\documentclass[12pt]{article}%revised version 2010-6-2
\def\q{\hfill\rule{1ex}{1ex}}
\def\0{\emptyset}

\def\QEDopen{{\hfill\setlength{\fboxsep}{0pt}\setlength{\fboxrule}{0.2pt}\fbox{\rule[0pt]{0pt}{1.3ex}\rule[0pt]{1.3ex}{0pt}}}}

\def\n{\noindent}

\usepackage{graphicx}
\usepackage{amsmath}
\usepackage{url}
\usepackage{graphicx}   %加图用
\usepackage{dsfont}
\usepackage{amsfonts}
\usepackage{amsmath}
\usepackage{bm}      %加粗
\usepackage{multirow,bigstrut,booktabs}    %画表格用
\usepackage{setspace}
\usepackage{indentfirst}
\usepackage{latexsym,  amsthm, graphicx, mathrsfs} %加载宏包
\usepackage{amssymb}
\usepackage{url} %引用相关？

\usepackage{geometry}
\begin{document}
\title{\bf Rainbow Independent Sets in Cycles
\thanks{
This work is partially supported by the National Natural Science Foundation of China (Grant 11771247 \& 11971158) and  Tsinghua University Initiative Scientific Research Program.}}
\author{{\small\bf Zequn Lv}\thanks{email:  lvzq19@mails.tsinghua.edu.cn
}\quad {\small\bf Mei
Lu}\thanks{email: lumei@tsinghua.edu.cn}\\
{\small Department of Mathematical Sciences, Tsinghua
University, Beijing 100084, China.}}

\date{}

\maketitle\baselineskip 16.3pt

\begin{abstract}

For a given class ${\cal C}$ of graphs and given integers $m \le n$,
let $f_{\cal C}(n,m)$ be the minimal number $k$ such that every $k$ independent
$n$-sets in any graph belonging to ${\cal C}$ have a (possibly partial)
rainbow independent $m$-set. In this paper, we consider the case ${\cal C}=\{C_{2s+1}\}$ and show that $f_{C_{2s+1}}(s, s) = s$. Our result is a special case of the conjecture (Conjecture 2.9) proposed by Aharoni et al in \cite{Aharoni}.

\end {abstract}

\vskip.3cm

\vskip.2cm

%{\bf Mathematics Subject Classification (2010):} 05C35, 05C38, 05C15

\vskip.2cm

\section{Introduction}
 Let $G=(V,E)$ be a finite, simple and undirected graph. A set of vertices in $G$ is called {\em independent} if it does not contain an edge of $G$. We call a set $S$ {\em $s$-set} if $|S|=s$. Now we cite the definition of {\em rainbow set} given in \cite{Aharoni}.

 {\bf Definition 1 \cite{Aharoni}} Let ${\cal F} =(F_1,\ldots , F_m)$ be a collection of (not necessarily distinct) sets. A {\em (partial) rainbow set} for ${\cal F}$ is  a set of the form
$R =\{x_{i_1} , x_{i_2},\ldots, x_{i_k}\}$, where $1 \le i_1 < i_2 < \ldots < i_k\le m$, and
$x_{i_j}\in F_{i_j}$ ($j\le k$). Here $R$ is a set, namely that the
elements $x_{i_j}$ are distinct.

When all sets in ${\cal F}$ are matchings of $G$, then $R$ is a (partial) rainbow matching for ${\cal F}$.
 There are some results about  rainbow matching when $G$ is a bipartite graph. Drisko \cite{Drisko} showed that every family of $2n- 1$ matchings of size $n$ in a bipartite graph where one side has size $n$ contains
a rainbow matching of size $n$. Aharoni and Berger \cite{Aharoni3} showed that
this is true for all bipartite graphs. A well-known conjecture given by Ryser \cite{Ryser} was
 that for odd $n$, every family of edge-disjoint $n$ perfect matchings of
$K_{n,n}$ has a rainbow matching of size $n$. In \cite{Aharoni2}, some results about  rainbow matchings for general graphs were given. %For the line graph H of a graph G, rainbow matchings in G correspond to
%rainbow independent sets in H.

Given
a family ${\cal F}$ of independent sets in a graph, a {\em rainbow independent set} is an
independent set whose vertices are contained in distinct independent sets in
${\cal F}$. For the line graph $H$ of a graph $G$, rainbow matchings in $G$ correspond to
rainbow independent sets in $ H$. So the results on rainbow matchings of a graph $G$ can be changed into rainbow independent sets of the  line
graph of $G$. In this paper, we consider  the rainbow independent sets.
 We first introduce parameters given in  \cite{Aharoni}. For a graph $G$ and integers $m \le n$, let $f_G(n,m)$ be the minimal
number $k$ such that every $k$ independent $n$-sets in $G$ have a partial
rainbow independent $m$-set. For a class ${\cal C}$ of graphs, let $f_{\cal C}(n,m) =\sup\{f_G(n,m) | G\in  {\cal C}\}$. This can be $\infty$. In the same paper, Aharoni, Briggs, Kim and Kim  \cite{Aharoni} established finiteness or infiniteness of $f_{\cal C}(n,m)$ for
many graph classes. In \cite{Kim}, Kim, Kim and  Kwon added two dense graph classes such that $f_{\cal C}(n,m)$ is finiteness.
In \cite{Aharoni}, Aharoni, Briggs, Kim and Kim  presented some conjectures. We will pay attention to  the case ${\cal C}=\{C_{2s+1}\}$ and consider the following conjecture (Conjecture 2.9 in \cite{Aharoni}).

\vskip.2cm
{\bf Conjecture 1 \cite{Aharoni}} {\em If $s < \frac{t}{2}$, then $f_{C_t}(s, s) = s$.}

\vskip.2cm
Note that $f_{C_t}(s, s) \ge s$. In this paper, we will show that Conjecture 1 holds if $t=2s+1$, that is, we have the following result.

{\bf Theorem 1 } {\em Let ${\cal I}=(I_1,I_2,\ldots,I_s)$ be a set of independent sets of size $s$ in $C_{2s+1}$ (not necessarily distinct). Then there is a rainbow independent set of size $s$ in $C_{2s+1}$, that is, $f_{C_{2s+1}}(s, s) = s$. }
\vskip.2cm

Let $V(C_{2s+1})=\Bbb Z_{2s+1}=\{1,2,\ldots,2s+1\}$ with $(i,i+1)\in E(C_{2s+1})$ (modulo $2s+1$) and ${\cal I}=(I_1,I_2,\ldots,I_s)$ be a set of independent sets of size $s$ in $C_{2s+1}$.  We define a bijection $f:\Bbb Z_{2s+1}\rightarrow \Bbb Z_{2s+1}$ such that $f(j)=2j-1$ for $1\le j\le 2s+1$ (modulo $2s+1$). Then $f^{-1}(I_i)$ consists of  $s$ consecutive vertices on $C_{2s+1}$ for all $1\le i\le s$. Thus Theorem 1 is  equivalent to the following theorem.

\vskip.2cm

{\bf Theorem 2 } {\em Let ${\cal F}=(V_1,V_2,\ldots,V_s)$, where $V_i$  consists of $s$ successive  vertices on the cycle $C_{2s+1}$ for $1\le i\le s$. Then there is a set of the form
$\{x_{i_1} , x_{i_2},\ldots, x_{i_s}\}$ such that
$x_{i_j}\in V_{j}$ ($1\le j\le s$) and $x_{i_j} x_{i_{j+1}}\in E(C_{2s+1})$ for all $1\le j\le s-1$.}

In Section 2, we will prove Theorem 2.

\vskip.2cm
\n{\large\bf 2.\ Proof of Theorem 2}
\vskip.2cm
Let $V(C_{2s+1})=\Bbb Z_{2s+1}=\{1,2,\ldots,2s+1\}$ with $(i,i+1)\in E(C_{2s+1})$ (modulo $2s+1$). Let $V_i=\{a_i,a_i+1,\ldots,a_i+s-1\}$  with $a_i\in  \Bbb Z_{2s+1}=V(C_{2s+1})$ for $1\le i\le s$. Note that if there is $i\in\{1,2,\ldots,s\}$ and $1\le j\le s-1$ such that $a_i+j>2s+1$, then $a_i+j$ is the vertex $(a_i+j)-(2s+1)$ on $C_{2s+1}$. Assume, without loss of generality, that $1= a_1\le a_2\le \ldots\le a_s\le 2s+1$. By Pigeonhole Principle, there is $k\in \{0,1,\ldots,2s\}$ such that $k\not\equiv a_i-1-i \pmod{2s+1}$ and $k\not\equiv a_i+s-i \pmod{2s+1}$ for all $1\le i\le s$. We have the following claim.

{\bf Claim 1} {\em (1) If $k>a_i-1-i$, then $k>a_j-1-j$ for all $1\le j\le i$.

(2) If $k>a_i+s-i$, then $k>a_j+s-j$ for all $1\le j\le i$.

(3) If $k<a_i-1-i$, then $k<a_j-1-j$ for all $i\le j\le s$.

(4) If $k<a_i+s-i$, then $k<a_j+s-j$ for all $i\le j\le s$.}

{\bf Proof of Claim 1} We just prove (1) and the other results can be proved by the same way.

Since $k>a_i-1-i$ and $ a_i\ge a_{i-1}$, we have $k>a_{i-1}-1-i$ which implies $k\ge a_{i-1}-1-(i-1)$.  Note that $k\not\equiv a_i-1-i \pmod{2s+1}$ for all $1\le i\le s$. Thus we have $k> a_{i-1}-1-(i-1)$. By repeating this process, we have (1) holds.\QEDopen

Now we complete the proof by considering the following two cases.

{\bf Case 1} $k<a_1+s-1$.

Since $a_1=1$, we have $k<s$. By Claim 1(4), we have $k<a_j+s-j$ for all $1\le j\le s$. Let $r=\max\{i~|~k>a_i-1-i\}$. Then $r\ge 1$. By Claim 1(3), we have $k<a_i-1-i$ for all $r+1\le i\le s$ if $r\le s-1$.

When $1\le i\le r$, we have $a_i-1<k+i<a_i+s$ which implies $k+i\in V_i$.

When $r+1\le i\le s$, we have $k+i-s<a_i+s-(2s+1)$ from $k<a_i-1-i$. We will prove that $a_i-1-(2s+1)<k+i-s$ by induction.

When $i=s$, we have $a_s-1-(2s+1)<k$ from $a_s\le 2s+1$ and $k\ge 0$. Suppose $a_i-1-(2s+1)<k+i-s$. That is $a_i+s-i\le k+(2s+1)$. Since $k\not\equiv a_i+s-i \pmod{2s+1}$ for all $1\le i\le s$, we have $a_i+s-i< k+(2s+1)$. Since $a_{i-1}\le a_i$, we have $a_{i-1}+s-i< k+(2s+1)$, that is, $a_{i-1}-1-(2s+1)< k+(i-1)-s$.

Now we have $a_i-1-(2s+1)<k+i-s<a_i+s-(2s+1)$ for all $r+1\le i\le s$. Then $k+i-s$ is a vertex in $V_i$  for all $r+1\le i\le s$. Thus $$\{k+i-s~|~r+1\le i\le s\}\cup\{k+i~|~1\le i\le r\}$$ consists of $s$ successive  vertices on the cycle $C_{2s+1}$.

{\bf Case 2} $k>a_1+s-1$.

In this case, $k>s$. Let $r=\max\{i~|~k>a_i+s-i\}$. Then $r\ge 1$ and $k<a_i+s-i$ for all $r+1\le i\le s$ if $r\le s-1$. By Claim 1(2),  $k>a_i+s-i$ for all $1\le i\le r$. Since $a_s\le 2s+1$ and $k>s$, we have $k>a_s-1-s$. By Claim 1(1), $k>a_i-1-i$ for all $1\le i\le s$.

When $r+1\le i\le s$, we have $a_i-1<k+i<a_i+s$ which implies $k+i\in V_i$.

When $1\le i\le r$, we have $k>a_i+s-i$ which implies $k+s+i>a_i-1+(2s+1)$. We will show that $k+s+i<a_i+s+(2s+1)$ by induction.

When $i=1$, we have $k+s+1<a_1+s+(2s+1)$ by $k<2s+1$. Suppose $k+s+i<a_i+s+(2s+1)$. Then $k-(2s+1)\le a_i-1-i$. Since $k\not\equiv a_i-1-i \pmod{2s+1}$ for all $1\le i\le s$, $k-(2s+1)< a_i-1-i$. Since $a_i\le a_{i+1}$, $k-(2s+1)< a_{i+1}-1-i$, that is, $k+s+(i+1)<a_{i+1}+s+(2s+1)$.

Now we have $a_i-1+(2s+1)<k+s+i<a_i+s+(2s+1)$ for all $1\le i\le r$. Then $k+s+i$ is a vertex in $V_i$  for all $1\le i\le r$. Thus $$\{k+i~|~r+1\le i\le s\}\cup\{k+s+i~|~1\le i\le r\}$$  consists of $s$ successive  vertices on the cycle $C_{2s+1}$. \q


\vskip.2cm
\begin{thebibliography}{99}

\bibitem{Aharoni} R. Aharoni, J. Briggs, J. Kim, M. Kim, Rainbow independent sets in certain classes of graphs, arXiv: 1909.13143v1.

\bibitem{Aharoni2} R. Aharoni, E. Berger, M. Chudnovsky, D. Howard, P. Seymour, Large rainbow
matchings in general graphs, European J. Combin.,
79(2019) 222-227.

\bibitem{Aharoni3} R. Aharoni and E. Berger, Rainbow matchings and matchings in $r$-partite $r$-graphs, Electron. J.
Combin. 16 (2009) \#R119.



\bibitem{Drisko} A.A. Drisko, Transversals in row-Latin rectangles. J. Combin. Theory
Ser. A, 84(2)(1998) 181-195.

%\bibitem{Erdos} P. Erd\v os, On sequences of integers no one of which dives the product of two others
%and some related problems, Mitt. Forsch.-Ins. Math. Mech. Univ. Tomsk 2 (1938),
%74-82.

%\bibitem{P1} P. Erd\v os, A. R\'enyi, On a problem in the theory of graphs. (Hungarian) Magyar Tud.
%Akad. Mat. Kutat$\acute{0}$ Int. K\"ozl. 7 (1962), 623-641.

%\bibitem{P2} P. Erd\v os, A. R\'enyi, V.T. S\'os, On a problem of graph theory, Studia Sci. Math.
%Hungar. 1 (1966), 215-235.

%\bibitem{Firke} F. Firke, P. Kosek, E. Nash, J. Williford, Extremal graphs without 4-cycles,
%Journal of Combinatorial
%Theory, Series B 103 (2013), 327-336.

%\bibitem{F} Z. F\"uredi, On the number of edges of quadrilateral-free graphs, Journal of Combinatorial
%Theory, Series B 68 (1996), 1-6.

\bibitem{Kim} J. Kim, M. Kim and O. Kwon, Rainbow independent sets on dense graph classes, arXiv: 2001.10566v2.





\bibitem{Ryser} H. J. Ryser, Neuere probleme der kombinatorik, Vortr\"age \"uber Kombinatorik
Oberwolfach, Mathematisches Forschungsinstitut Oberwolfach,
July 1967.

\end{thebibliography}
\end{document}